\documentclass{amsart}
\usepackage{amsmath}
\usepackage{amssymb}
\usepackage{amsmath,amsthm,amssymb,amsfonts}

\usepackage{framed}

\usepackage{hyperref}

\usepackage{graphicx}

\usepackage{url}

\newtheorem{Def}{Definition}
\newtheorem{Th}{Theorem}
\newtheorem{Lm}{Lemma}
\newtheorem{Corol}{Corollary}
\newtheorem{Prop}{Proposition}
\newtheorem{Exa}{Example}
\newtheorem{Con}{Conjecture}

\def\cl {\mbox{\rm cl}\,}

\def\R{{\mathbb{R}}}

\def\proof{{\bf Proof.\ }}

\def\ep{\varepsilon}

\def\L {\mathcal L}
\def\d {\partial}
\def\i {\mbox{\rm int}\,}
\def\ll {\lambda}

\def\x {\overline{x}}

\begin{document}
\title[Gradient methods]{Gradient projection  and  conditional gradient methods for constrained nonconvex minimization
}

\author{ M. V. Balashov, B. T. Polyak, A. A. Tremba} 
\address{V. A. Trapeznikov Institute of Control Sciences of Russian Academy of Sciences,
65 Profsoyuznaya street, Moscow 117997, Russia.
}
\email{balashov73@mail.ru, boris@ipu.ru, atremba@ipu.ru}


\maketitle              

\def\R {\mathbb{R}}

\begin{abstract}
    Minimization of a smooth function on a sphere or, more generally,
    on a smooth manifold, is the simplest non-convex optimization
    problem. It has a lot of applications. Our goal is to propose a
    version of the gradient projection algorithm for its solution and
    to obtain results that guarantee convergence of the algorithm
    under some minimal natural assumptions. We use the
    Le\v{z}anski-Polyak-Lojasiewicz condition on a manifold to prove
    the global linear convergence of the algorithm.
    Another method well fitted for the problem is  the conditional gradient (Frank-Wolfe) algorithm.
    We examine some conditions which guarantee global convergence of full-step version of the method
    with linear rate.

\end{abstract}

\section{Introduction}

Consider  minimization  of a smooth function $f(x)$
on a closed set $Q$ in the Euclidean space $(\R^n,\|\cdot\|)$

\begin{equation}\label{sphere}
\min_{x\in Q} f(x).
\end{equation}

Traditionally the set $Q$ and the function $f$ are assumed to be convex; in such convex setting the problem is well studied and numerous algorithms are known, see e.g.
\cite{Bertsekas,Boyd, Nesterov1,NesterovNemir} for details.
We plan to address the situation when the function or/and the set are nonconvex.

The function $f$ under consideration is smooth with the Lipschitz continuous gradient, but nonconvex.
Regarding the set $Q$,  we mostly suppose that it is proximally smooth \cite{Vial, Clarke, Clarke1}.
In particular we consider the next important cases:\\
1) $Q=S_R=\{x \in \mathbb{R}^{n} \ |\ ||x||=R\}$ (minimization on the sphere),\\
2) $Q = \{ x\in\R^n\ |\ g(x)=0\}$, $g:\R^n\to\R^m$ (equality type constraints),\\
3) $Q$ is the boundary of a strongly convex set $B\subset\R^n$.\\

Minimization on the sphere has numerous applications, for
instance finding minimal eigenvalue of a symmetric matrix
$A\in\R^{n\times n}$ (then $f(x)=(Ax,x)$, $Q=S_1$) or choosing
step-size in trust-region methods \cite{trust-region}. The set
$Q$ is obviously nonconvex thus (\ref{sphere}) is an example of
nonconvex optimization problems. The pioneering work in the field
is \cite{L}. Special case of  problem (\ref{sphere}) (for
quadratic $f(x)$) has been studied by Hager \cite{Hager}; the
solution can be reduced to solving 1D equations. Later
publications include \cite{Absil-2008,
Neto,Lacoste-2016,Udriste,Weber-Sra-}, in most of them $Q$ is a
smooth (Riemannian) manifold. But in general, research in
nonconvex optimization is much less intensive than in convex
case. The main approaches use generalized convexity on the set
and/or consider geodesic-related steps. In contrast, we use
neither of these. There are numerous methods for optimization
with equality type constrains, see e.g \cite[Chapter
7]{Polyak_book}, \cite{Bertsekas}. However, most of them generate
points which are not admissible ($x_k\notin Q$) while our purpose
is to develop methods with admissible iterations.

The contribution of the present paper is triple:
\begin{enumerate}
    \item We propose a new approach to the gradient projection algorithm for constrained optimization,
    based on the idea of upper approximation of the objective function. The resulting method is
    the gradient projection algorithm with constant step-size; it differs from versions proposed in \cite{L}.
    Moreover we prove its convergence to a stationary point without any convexity-like assumptions.
    \item We generalize well known in unconstrained optimization the Polyak-Loja\-sie\-wicz condition on
    the class of problems (\ref{sphere}). Under this assumption we prove
    linear convergence of the gradient projection algorithm (with projection on the 
    tangent subspace and a variant combined with the Newton method) to a global extremum under assumption of proximal smoothness of the manifold $Q$.
    As an example we consider a quadratic form on the unit sphere.
    \item For \emph{approximately linear} objective functions we propose a new version of the Frank-Wolfe
    (conditional gradient) method and establish its linear convergence to a global minimum in problem (\ref{sphere}). We prove
    linear convergence of the method for a surface, which is the boundary of a strongly convex set. Note, that such surface is not
    necessary smooth.
\end{enumerate}

In the paper \cite{Balashov2017} linear convergence of the gradient projection algorithm was
proved for a proximally smooth set with constant $R$ and for a strongly convex function with constant of strong convexity $\varkappa$ and Lipschitz constant $L$
under assumption $\frac{L}{\varkappa}<R$. The last inequality is essential for linear convergence
of the method.
In subsection \ref{SS-AP} we prove convergence of the standard gradient projection algorithm
to a stationary point of the problem (\ref{sphere}) for any function with the Lipschitz
continuous gradient and for any proximally smooth set.

In subsection \ref{11} we extend  the Polyak-Lojasiewicz condition (well known in
unconstrained minimization)
\cite{Lojasiewicz,Polyak-63,Karimi} for the constrained case
(\ref{sphere}) with differentiable function and $C^1$ smooth
manifold $Q$. We want to pay attention that this
generalization is in fact some variant of the error bound
condition for the case of smooth set $Q$.

On the base of the new definition from subsection \ref{11} we
prove in subsection \ref{GPA-Sp} convergence of the gradient projection algorithm for
$Q=S_1$. For the case of general $C^1$ smooth
and proximally smooth manifold of codimension 1 we prove it in subsection
\ref{GPA-Any}. In contrast with other approaches our algorithms
represent variants of the gradient projection algorithm with
admissible points $x_k\in Q$ and linear rate of convergence.

In subsection \ref{GPA-NM} we consider the situation when the gradient projection
algorithm can be finalized by use of the Newton method. This is standard practice,
except that we are dealing with nonconvex problem.

Concluding section \ref{GPA-sec}, we prove the
Le\v{z}anski-Polyak-Lojasiewicz condition for the quadratic function on the sphere.
Thus we extend the result \cite{gao-etal-exponent-2016}
by clarifying  constant $\mu$ in this condition. It is essential for estimate of the error
for the gradient projection algorithm in the case $Q=S_1$.

In section \ref{FW-sec} we consider application of the Frank-Wolfe method
for solving the problem (\ref{sphere}) in the case when $Q$ is the boundary
of a strongly convex set of radius $R$.
It is well known that if the set $Q$ is a convex compact and $f$
is a Lipschitz differentiable convex function then the method
converges (with respect to the objective function) with sublinear
rate. In \cite{FW_1} the authors discuss the choice of step-size
in the method. In \cite{Lacoste-2016} the author proved that under
certain assumptions
in the case of convex compact set $Q$ and for a (nonconvex) function with the Lipschitz continuous gradient the Frank-Wolfe method converges to a stationary point in the problem
(\ref{sphere}) with sublinear rate.

We prove linear rate of convergence with respect to the point. In subsection
\ref{FW-ssec1} we proved it for so-called linear approximative function and in subsection
\ref{FW-ssec2} we proved it for a function with gradient domination.
In fact in both subsections the general idea consists in the fact that radius of ''curvature''
of the level sets for our function is larger than the radius of strong convexity. This leads
to the results. We only want to point out that we take the notation of \it radius of ''curvature'' \rm
in the sense of supporting principles for proximally smooth and strongly convex sets.

For completeness in Appendix we prove necessary condition of
minimum in the problem (\ref{sphere}) for a proximally smooth
set $Q$ and a function with the Lipschitz continuous gradient.

All mentioned parts are gathered together by the possibility of certain
``spherical" approximation (via the supporting principles) for the set/surface $Q$
in problem (\ref{sphere}).\\

The most part of mentioned results takes place in the case of the real Hilbert space.
Sometimes obvious patches should be applied in the infinite dimension case,
for example compactness of the set $Q$ in Theorem \ref{T1-1}.

\section{Definitions and main notations}

Let $B_R(a) = \{ x\in\R^n\ |\ \| x-a\|\le R\}$ be the ball with center $a\in\R^n$
and radius $R>0$.
 For a set $Q\subset \R^n$ the sets $\cl Q$, $\i Q$, $\d Q$ \it are the closure, the interior and
the boundary \rm of $Q$, respectively. We also denote by $\d Q$ the edge of the surface $Q$. 

Let $P_Q$ be \textit{the operator of metric projection} on the set $Q$, i.e. $P_Q(x)=\{y\in Q: ||y-x||=\inf_{z\in Q} ||x-z|| \}$. In general, $P_Q$ can be set-valued for nonconvex sets, but for proximally smooth sets (see below) it is single-valued (provided $x$ is close enough to $Q$).

For a closed set $Q\subset\R^n$ \it the normal cone of proximal normals
\rm (or simply the normal cone) at a point $x\in Q$ is defined as follows \cite{Clarke1}
$$
N(Q,x) = \{ p\in\R^n\ |\ \exists \delta>0\ P_Q(x+\delta p)= x \}.
$$
If the set $Q$ is convex, then $N(Q,x)$ coincides with normal cone
in the sense of convex analysis.

A closed set $Q\subset \R^n$ is called \it proximally smooth \rm
with constant $R>0$ \cite{Vial, Clarke, Clarke1} if the distance
function $\varrho (x,Q) =\varrho _Q(x)=\inf\limits_{a\in Q}\|x-a\|$
is continuously Frechet
differentiable on the set $U_Q(R) = \{ x\in\R^n\ |\ 0<\varrho
_Q(x)<R\}$. The equivalent properties for a proximally smooth set
with constant $R$ are

1) $P_Q:U_Q(R)\to Q$ is a single-valued continuous function,

2) \it supporting principle: \rm $p\!\in\! N(Q,x)$, $\|p\|\!=\!1$, if and only if $Q\cap\i
B_R(x+Rp)=\emptyset$.

Note that the mapping $Q\ni x\to N(Q,x)$ is upper semicontinuous for a proximally smooth set $Q$
with constant $R$. For a point $x\in\i Q$ we have $N(Q,x)= 0$ and it is sufficient to prove upper semicontinuity on the boundary $\d Q$.
Choose $x_k\in\d Q$, $x_k\to x\in \d Q$, and $p_k\in N(Q,x_k)$, $\| p_k\|=1$, with
$p_k\to p$. Suppose that $p\notin N(Q,x)$, the last means that $N(Q,\cdot)$ is not
upper semicontinuous at the point $x$. By the supporting principle for proximally smooth sets
$$
\exists\, \x \in Q\cap \i B_R(x+Rp).
$$
The sets $\i B_R(x_k+Rp_k)$ converge to the set $\i B_R(x+Rp)$ in the Hausdorff metric,
thus $\x \in Q\cap \i B_R(x_k+Rp_k)$ for sufficiently large $k$. A contradiction.

If we consider a continuously differentiable $(n-1)$-dimensional surface $Q$ without edge which is also proximally smooth with constant $R>0$, then for any point $x\in Q$ the normal cone $N(Q,x)$ is
1-dimensional subspace. If $\| p\|=1$ and $p\in N(Q,x)$ then the surface $Q$ is trapped between the supporting
spheres (see supporting principle 2)):
\begin{equation}\label{BetweenS}
Q\subset \left( \R^n\backslash \i B_R(x-Rp)\right)\cap \left( \R^n\backslash\i B_R(x+Rp)\right).
\end{equation}

A closed convex set $B\subset \R^n$ is called \it strongly convex
of radius $r>0$ \rm if it can be represented in the form
$B=\bigcap\limits_{x\in X}B_r(x)$ \cite{Vial,PBMsb}.
There are few equivalent properties for strong convexity:\\
1) A convex
compact set $B\subset\R^n$ is strongly convex of radius $r$ if and only if
for any pair of points $x,y\in B$ the ball with center
$\frac12 (x+y)$ of radius $\frac{1}{8r}\| x-y\|^2$ belongs to $B$
\cite{PolykLevitin, Polyak66,Balashov12}. \\
2) Another equivalent property for
strong convexity is \it supporting principle: \rm for any $x\in \d
B$ and $p\in N(B,x)$, $\| p\|=1$, we have
$$
B\subset B_r(x-rp).
$$
3) The set $B\subset\R^n$ is
strongly convex of radius $r$ if and only if for any unit vectors
$p,q$ and for supporting elements $x_p=\arg\max\limits_{x\in
B}(p,x)$, $x_q=\arg\max\limits_{x\in B}(q,x)$ we have the next
inequality \cite[Proposition 2.8]{Vial}
$$
\| x_p-x_q\|\le r\| p-q\|.
$$

If $Q$ is a smooth manifold, then the subspace $T_x$ is \textit{the tangent subspace} to  $Q$ at a point $x\in Q$, i.e.
$$
T_x = \{ v\in\R^n:  \varrho (x+tv, Q)=o_{v}(t) \},
$$
where $\lim\limits_{t\to+0}\frac{o_{v}(t)}{t}=0$ for all $v\in
T_x$.

Define $\L_f(\beta) = \{ x\in\R^n\ |\ f(x)\le \beta\}$ \it the lower level set
\rm of the function $f$.

We  say that the function $f$ has the \it  Lipschitz continuous gradient
\rm with constant $L_1>0$ if
\begin{equation}\label{Lip}
\|f'(x)-f'(y)\|\le L_1\|x-y\|, \quad \forall x,y.
\end{equation}

It is well known that (\ref{Lip}) implies the upper bound for $f$
\cite[Lemma 1.2.3]{Nesterov1}:
\begin{equation}\label{key}
\left|f(y) - f(x) - (f'(x),y-x)\right|\le\frac{C}{2}||y-x||^2, \quad \forall x,y,
\end{equation}
for any $C\ge L_1$.

We write $f(t)\asymp g(t), t\in \R^1$, $t\to +0$, if there exist $0<C_1<C_2$ and $\delta>0$
such that $C_1 f(t)\le g(t)\le C_2f(t)$ for all $t\in (0,\delta)$.

For a differentiable vector function $g(x) = (g_1(x),\dots, g_m(x))^T$ denote the Jacobi
matrix as
$$
g'(x) = (g_1'(x)\dots g_m'(x))
$$
and we treat $g_{i}'(x)$ as columns.

The set $\Omega \subset Q$ is the set of \it
stationary points \rm of the differentiable function $f$ on the
set $Q$ (which is associated with the problem (\ref{sphere})) if
for any point $x\in \Omega$ we have $f'(x)\in - N(Q,x)$. The last
inclusion is necessary condition of optimality for a proximally
smooth set $Q$ and a smooth function $f$. We prove it in Appendix for
completeness.

\section{The gradient projection algorithm}\label{GPA-sec}

The gradient projection algorithm
for \eqref{sphere} in convex case has been proposed in
\cite{Goldstein,PolykLevitin}. The simplest version (with constant step-size) looks as follows:
for an iteration $x_k\in Q$ it generates the new point $x_{k+1}$ as the minimizer
of the upper bound on $Q$ (with $C_1\ge L_1$) or, equivalently, as projection of the gradient step on $Q$:
\begin{framed}
\textbf{Gradient Projection Algorithm (GPA1)}
\vspace{3mm}

\textbf{Step 1.}
Choose a constant $C_1 > 0$, initial point $x_0 \in Q$ and put $k = 0$.

\textbf{Step 2.}
Repeat
\begin{equation}\label{upper}
\begin{array}{rl}
x_{k+1} \!\!\!\!& = \displaystyle \arg\min_{x\in Q} \left\{f(x_k) + (f'(x_k), x-x_k) + \frac{C_1}{2}||x-x_k||^2\right\} = \\[3mm]
& = \displaystyle P_Q\left( x_k-\frac{1}{C_1}f'(x_k)\right),
\end{array}
\end{equation}
\end{framed}
The condition $C_1\ge L_1$ is equivalent to choice of the constant step-size $t \leq \frac{1}{L_1}$.
Below we shall consider the extensions of the method for nonconvex set $Q$ and nonconvex function $f(x)$.

One of the possible ways for extension
 is \emph{gradient projection along geodesics}  proposed by Luenberger \cite{L}. It is not hard to design geodesics on the sphere (arcs of big circles on the sphere),
but the original algorithm in \cite{L} requires one-dimensional minimization on each iteration.
Another problem is that in the case of an arbitrary manifold $Q$ construction of geodesics is a hard procedure.
Thus we avoid geodesics and try to deal with gradient projection method with constant step-size in the form \eqref{upper}.

\subsection{The case of an arbitrary proximally smooth set $Q$. General algorithm}\label{SS-AP}

The next result shows that for any function $f:\R^n\to\R$ with the Lipschitz continuous
gradient and for any proximally smooth set $Q\subset\R^n$  iterations of the standard gradient
projection algorithm \eqref{upper} are well-defined and
 converge to a stationary point of the function $f$ on the set $Q$
for the appropriate choice of the step-size.

\begin{Th}\label{T1-1} \it
    Let $Q\subset \R^n$ be a bounded proximally smooth set with constant $R>0$.
    Suppose that the function $f:\R^n\to\R$ is the Lipschitz continuous with constant $L>0$
    and its gradient $f'$ is also the Lipschitz continuous with constant $L_1>0$.
    Take $C>0$ with $\frac{L}{C+L_1}<R$. Then
    for any $x_0\in Q$
    GPA1 with $C_1 = C + L_1$
converges to the set of stationary points $\Omega$: $\lim\limits_{k\to \infty}\varrho_{\Omega}(x_k) = 0$ and
    $$
    f(x_{k+1})+\frac{C}{2}\| x_{k+1}-x_k\|^2\le f(x_k)
    $$
    for all $k\ge 0$.\rm
\end{Th}
\proof Define for each natural $k$ the function
$$
\psi_k(x) = f(x_k)+
(f'(x_k),x-x_k)+\frac{C+L_1}{2}\| x-x_k\|^2.
$$
It's easy to see (due to the Lipschitz continuity of gradient $f'$) that
$$
\psi_k(x) \ge f(x)+\frac{C}{2}\| x-x_k\|^2
$$
for all $x\in\R^n$ and
$$
\varrho_Q\left(x_k-\frac{1}{C+L_1}f'(x_k) \right)\le \left\| \frac{1}{C+L_1}f'(x_k)\right\| \le \frac{L}{C+L_1}<R.
$$
Hence the distance from the point $x_k-\frac{1}{C+L_1}f'(x_k)$ to the set $Q$ is less than
$R$ and the metric projection $x_{k+1}$ is defined uniquely by the definition of proximally smooth set.
We have
$$
f(x_k)=\psi_k (x_k)\ge \psi_k (x_{k+1})\ge f(x_{k+1})+\frac{C}{2}\| x_{k+1}-x_k\|^2.
$$

Assume that $\lim_{k\to\infty}\varrho_{\Omega}(x_k)\ne 0$ for some
sequence $\{x_k\}\subset Q$ which is generated by the gradient
projection algorithm.
Then
there is a number $\ep>0$
and a subsequence $\{x_{k_m}\}\subset \{x_k\}$
with $\varrho_{\Omega}(x_{k_m})\ge\ep$ for all $m$.
Consider a
converging subsequence of the sequence $\{x_{k_m}\}$
(that again is denoted by $\{x_{k_m}\}$)
and $x_* = \lim\limits_{m\to\infty}x_{k_m}$.
Then from the necessary conditions of minimum of the function $\psi_{k_m}$ on the set $Q$ we get
$$
\psi'_k(x_{k_m +1 })\in - N(Q,x_{k_m +1 }),
$$
in other words
$$
f'(x_{k_m})+(C+L_1)(x_{k_m +1 }-x_{k_m})\in - N(Q,x_{k_m +1 }).
$$
Passing to the limit as $m\to \infty$, using upper semicontinuity of the normal cone $N(Q,\cdot)$
and the property $\lim_{m\to \infty}(x_{k_m +1 }-x_{k_m})=0$
we have
$$
f'(x_*)\in - N(Q,x_{*}).
$$
Thus $\lim\limits_{m\to\infty} \varrho_{\Omega}(x_{k_m})=0$, a contradiction.
\qed

With the help of Theorem \ref{T1-1} we can find a stationary point with error $\ep>0$,
namely we can find a point $x\in Q$ with $\varrho (0, -f'(x)-N(Q,x))<\ep$.

\begin{framed}
\textbf{Stationary-point Algorithm}

\textbf{Step 1.}
Choose $\ep > 0$ and $C > \max\{0, \frac{L}{R}-L_1 \}$. Put $\Delta f= \sup_Q f - \inf_Q f$, $\delta=\frac{\ep}{C+2L_1}$,
and set $k = 0$.

\textbf{Step 2.} Perform a Step~2 \eqref{upper} of GPA1 with $C_1 = C + L_1$.

\textbf{Step 3.} If $\|x_k - x_{k+1}\| \geq \delta$, increase $k$ and continue to the Step~2.\\
Otherwise stop the algorithm and return $x_{k+1}$ as the solution.

\end{framed}

The algorithm do at most $N = \big\lfloor \frac{2 \Delta f}{C\delta^2} \big\rfloor + 1$ steps.
Assume the contrary, that $\| x_{k+1}-x_k\| \geq \delta$ for $k=0,\dots,N-1$.
Then assumptions of Theorem~\ref{T1-1} holds, thus
$\frac{C}{2}\delta^2\le f(x_k)-f(x_{k+1})$ and after $N$ steps we get
$$
f(x_0) - f(x_N) \geq N\frac{C}{2}\delta^2 > \Delta f,
$$
i.e. a contradiction.
When the algorithm stops with $\|x_k - x_{k+1}\| < \delta$, by the optimality condition for the function $\psi_k$ we get
$$
f'(x_{k})+(C+L_1)(x_{k+1}-x_{k})\in - N(Q,x_{k+1}).
$$
Using the Lipschitz continuity of $f'$ we obtain
$$
\varrho (f'(x_{k+1}),-N(Q,x_{k+1}))\le \delta (C + 2L_1)=\ep.
$$

Note that the parameter $C$ of the
step $t=\frac{1}{C+L_1}$ depends on Lipschitz constant $L$ of the function
$f$ because the point $x_k-t f'(x_k)$
should not go very far from the set $Q$.
Moreover, conditions on constant $C$ mean that the step-size $t$ satisfies the inequality $t<\min\left\{\frac{R}{L},\frac{1}{L_1}\right\}$.

\begin{Exa}\label{Sphere1}
Sometimes projection can be found explicitly, for example for the unit sphere
$Q=S_1$:
\begin{equation}\label{alg-with-L}
x_{k+1} = \frac{ x_k -t f'(x_k)}{\| x_k - t f'(x_k)\|},
\end{equation}
thus the algorithm is the gradient-projection method with
constant step-size.
\end{Exa}

By Proposition \ref{NecMin} the first-order optimality condition in (\ref{sphere}) for $Q=S_1$
means $f'(x)=\lambda x, \, \|x\|=1, \lambda\in \mathbb{R}$, and it can be immediately rewritten in the form
\begin{equation}\label{opt1}
\|(I - x x^T) f'(x)\|=0,
\end{equation}
(here $I$ is the unit matrix) or
as
\begin{equation}\label{opt}
\|L_1 x - f'(x)\| = |( L_1 x - f'(x), x)|.
\end{equation}

\subsection{The Le\v{z}anski-Polyak-Lojasiewicz condition on a manifold.}\label{11}

Now the main task is to propose conditions which guarantee convergence of the method to the global minimum in the problem
(\ref{sphere}) and to estimate the rate of convergence.

In unconstrained minimization we have such powerful tool as convexity; gradient method for convex functions converges to global minima while for strongly convex functions one has linear rate of convergence. There are extensions of convexity for minimization on manifolds, see e.g. the monograph
\cite{Udriste} and the paper \cite{Neto}. Unfortunately there exist no (globally) convex functions on compact manifolds \cite{Neto}, thus we need  some other tools.

However in the unconstrained case there are conditions which validate convergence for nonconvex functions. Probably the first one is due to T. Le\v{z}anski \cite{Lez0,Lez}. He considered
a problem of unconstrained minimization for a Lipschitz differentiable function $f$ such that there exists
a positive continuous function $\varphi$ with
$$
\| f'(x)\|\ge \varphi (f(x)-f_0),\quad \forall x\in\R^n,\qquad f_0=\min\limits_{x\in\R^n}f(x),\quad
\int\limits_{0}^{f(x)-f_0}\frac{ds}{\varphi (s)}<\infty.
$$

Under these assumptions he proved for $\varphi (s)=c\sqrt{s}$ the convergence of the gradient descent
algorithm
$
x_{k+1}=x_k-t_kf'(x_k)
$
with linear rate.
The same assumption
$$\|f'(x)\|^2\ge \mu (f(x)-f_0)$$
(where $\mu>0, f_0=\min f(x)$) was considered in
\cite{Lojasiewicz,Polyak-63}. Sometimes \cite{Karimi} this
 is referred as the Polyak-Lojasiewicz condition (works of Le\v{z}anski were not widely known). Thus it is fair to call the above condition as Le\v{z}anski-Polyak-Lojasiewicz (LPL) one.

Analogously we can propose the analog of LPL condition for the constrained minimization of a differentiable
function $f:\R^n\to\R$ on a smooth
manifold $Q\subset\R^n$. Define by $P_{T_x}$ the metric projection on the tangent subspace
$T_x$ to the manifold $Q$ at the point $x\in Q$.
Note that $N(Q,x)$ is the polar cone (subspace) for the tangent space $T_x$.

\begin{Def}\label{PL-general}
Let $Q\subset\R^n$ be a $C^1$ manifold and $f:\R^n\to\R$ be a differentiable function. Let $\mu>0$, $\alpha\ge 1$, $\beta\in\R^n$,
$f_0 = \min_Q f$. We shall say that the Le\v{z}ans\-ki-Po\-lyak-Lo\-ja\-sie\-wicz LPL condition with exponent $\alpha$ holds
for the function $f$ on the set $Q$ if
\begin{equation}\label{PLG}
\| P_{T_x}f'(x)\|^\alpha\ge \mu (f(x)-f_0)
\end{equation}
for all $x\in Q\cap \L_f(\beta)$.
\end{Def}

If $\alpha=2$ then we shall call (\ref{PLG}) simply the \textit{LPL condition} for the function $f$
on the manifold $Q$.

We want to admit that we consider such manifold $Q$ that
its edge $\d Q$ has empty intersection with the set $\L_f(\beta)$.
For example, $Q$ can be a manifold without edge.

 Note that if $Q$ is given by the system $S=\{ x\in\R^n\ |\ g(x)=0\}$ of
full rank then $P_{T_{x}} = I - g'(x)(g'(x)^Tg'(x))^{-1}g'(x)^T$ for all $x\in Q$.
Here $I$ is the identity operator in $\R^n$.

\begin{Exa}\label{LPLExa}
 In the case $Q=S_1$ LPL condition reads
\begin{equation}\label{PL}
\|(I - x x^T) f'(x) \|^2 \ge \mu (f(x)-f_0),
\end{equation}
 for all $x\in S_1\cap\L_f(\beta)$.\\
\end{Exa}
Later we shall consider quadratic case ($f(x)$ is a quadratic function) and confirm fulfillment of condition \eqref{PL}.

Now consider a special 2D example which exhibits possible situations.
\begin{Exa}\label{LPL2D}
 Let $Q=\{ (x,y)\in\R^{2}\ |\
 x^{2}+(y-\frac12)^{2}=\frac14\}$  (the set $Q$ is the circle with
 center $(0,\frac12)$ and radius $r=\frac12$). Let $f(x,y) =
 y-px^{2}$, where $p\in (0,1]$ is a parameter. We have
 $$
 (0,0) = \arg\min\limits_{(x,y)\in Q}f(x,y)
 $$
 and $f_{0} = f(0,0)=0$.
\end{Exa}

 3.1. Suppose that $p\in (0,1)$. Consider $t>0$ and
 $a_{0}=(x_{0},y_{0})\in Q$, $x_{0}>0$, with $f(x_{0},y_{0}) =
 y_{0}-px_{0}^{2}=t$. Note that $x_{0}\to 0$ is equivalent to
 $t\to 0$. The angle between tangent lines to the circle $Q$ and
 curve $f(x,y)=y-px^{2}=t$ at the point $a_{0}$ asymptotically
 equals $2(1-p)x_{0}$ when $x_{0}\to 0$. Substituting
 $y_{0}=px_{0}^{2}+t$ to the equation
 $x_{0}^{2}+(y_{0}-\frac12)^{2}=\frac14$ we obtain that
 $x_{0}^{2}\asymp t$, $x_{0}\to 0$. Hence $\|
 P_{T_{a_{0}}}f'(a_{0})\|\asymp x_{0}$, $x_{0}\to 0$. From the
 other hand $f(a_{0})-f_{0}=t\asymp x_{0}^{2}$. Thus the exponent
 in the LPL condition equals $\alpha = 2$.

 3.2. Suppose that $p=1$. Consider $t>0$ and $a_{0}=(x_{0},y_{0})\in
 Q$, $x_{0}>0$, with $f(x_{0},y_{0}) = y_{0}-x_{0}^{2}=t$. The
 angle between tangent lines to the circle $Q$ and curve
 $f(x,y)=y-x^{2}=t$ at the point $a_{0}$ asymptotically equals
 $x_{0}^{3}$ when $x_{0}\to 0$ and $x_{0}^{4}\asymp t$ when
 $x_{0}\to 0$. Hence $\| P_{T_{a_{0}}}f'(a_{0})\|\asymp
 x^{3}_{0}$, $x_{0}\to 0$. From the other hand
 $f(a_{0})-f_{0}=t\asymp x_{0}^{4}$ when $x_{0}\to 0$. Thus the
 exponent in the LPL condition equals $\alpha = \frac43$.

 We conjecture that these two situations are the only possible ones for quadratic objects.

\begin{Con} Let $A_{i}\in \R^{n\times n}$ be a symmetric matrix, $b_{i}\in \R^{n}$,
$i=1,2$ and $c\in\R$. Suppose that the set $Q$ is a quadric nonempty
surface, i.e. $Q=\{ x\in \R^{n}\ |\ (x,A_{1}x)+(b_{1},x)+c=0\}$,
$f(x)=(x,A_{2}x)+(b_{2},x)$ and there exists a unique point of the
global minimum for the problem $\min_{Q}f$. Then the exponent
$\alpha$ in the LPL condition near the global
minimum equals $2$ or $\frac43$.
\end{Con}

\subsection{The gradient projection algorithm on the unit sphere}\label{GPA-Sp}

Next we consider a special case when $Q=S_1$. In this case all projections can
be explicitly calculated.

\begin{Th}\label{T2} Suppose that $f$ is the Lipschitz function with constant $L>0$,
$f'$ is the Lipschitz function with constant $L_1>0$ and $0<\mu<2( L_1+L)$.
Under LPL condition (\ref{PL}) algorithm \eqref{alg-with-L} with
$t=\frac{1}{L_1}$ converges to a point of minimum $x_*$ with
linear rate.
\end{Th}

\proof
First, let's describe quantitative connection between the
term $\|(I - x x^T) f'(x)\|$ of optimality condition \eqref{opt1} and the residual
$z = z(x)= \|L_1 x - f'(x)\| - (L_1 x - f'(x), x) \geq 0$.

After simple arithmetical calculations we get
$$
z \cdot (\|L_1 x - f'(x)\| + (L_1 x - f'(x), x)) = \|(I - x x^T) f'(x)\|^2,
$$
and
$$
z = \frac{\|(I - x x^T) f'(x)\|^2}{\|L_1 x - f'(x)\| + (L_1 x - f'(x), x)}
\geq \frac{\|(I - x x^T) f'(x)\|^2}{2\|L_1 x - f'(x)\|},
$$
and calculations are well defined for any nonstationary point $x$, see conditions \eqref{opt1}, \eqref{opt}.

Fix a point $x_0\in S_1\cap\L_f(\beta)$.
We have $\| x_{k+1}-x_{k}\|^2=$
$$
=\left\|\frac{ L_1 x_k - f'(x_k)}{\| L_1 x_k - f'(x_k)\|} -
x_k\right\|^2 = 2 - 2 \frac{ (L_1 x_k - f'(x_k), x_k)}{\| L_1 x_k
- f'(x_k)\|} = 2 \frac{z_k}{\| L_1 x_k - f'(x_k)\|}.
$$

By the previous formula, (\ref{key}) and definition of $x_{k+1}$
by \eqref{alg-with-L} with $t=\frac{1}{L_{1}}$ the next estimate
holds
$$
f(x_{k+1})-f(x_{k})\le
(f'(x_{k}),x_{k+1}-x_{k})+\frac{L_{1}}{2}\| x_{k+1}-x_{k}\|^{2} =
-z(x_k).
$$
Hence
$$
f(x_{k+1}) - f(x_k) \leq -z(x_k) \leq -\frac{\|(I - x_k x_k^T) f'(x_k)\|^2}{2\|L_1 x_k - f'(x_k)\|}
$$
and thus $f(x_k)\le \beta$ implies $f(x_{k+1})\le f(x_k)\le \beta$.
Denoting $\varphi_k = f(x_k) - f_0$ and using condition
\eqref{PL} we get
$$
\varphi_{k+1} - \varphi_k \leq -\frac{\mu}{2\|L_1 x_k - f'(x_k)\|} \varphi_k,
$$
or
\begin{equation} \label{formula-for-q}
\varphi_{k+1} \leq \Big(1 - \frac{\mu}{2\|L_1 x_k - f'(x_k)\|}\Big) \varphi_k,
\end{equation}
From the latter inequality follows
$$
\varphi_{k+1} \leq q^k \varphi_0,
$$
where $q = 1 - \frac{\mu}{2(L_1 + \max_{x \in S_1} \|f'(x)\|)} \le
1 - \frac{\mu}{2(L_1 + L)}\in (0,1)$.

Now prove the convergence with respect to $x$. Note that obvious condition $q\ge 0$ implies $\|L_1 x_k - f'(x_k)\|\ge \mu / 2$. Thus for $\delta_k=x_{k+1}-x_k$ and $\|x_k\| = \|x_{k+1}\| = 1$ we get
\begin{gather*}
\|\delta_k\|^2 = 2 \frac{z_k}{\| L_1 x_k - f'(x_k)\|} \leq
\frac{4}{\mu} z_k \leq \frac{4}{\mu}(f(x_{k}) - f(x_{k+1}))
\leq\frac{4}{\mu} \varphi_{k} \leq \frac{4\varphi_0}{\mu} q^k.
\end{gather*}
Hence $\|\delta_k\|\leq \sqrt{\frac{4\varphi_0}{\mu}}\, q^{\frac{k}{2}}$ and $\{x_k\}$ is the Cauchy sequence. This  implies its convergence to a point $x_*$ with linear rate, while inequality $\varphi_{k+1} \leq q^k \varphi_0$ and continuity of $f(x)$ provides $f(x_*)=f_0$. \qed

\subsection{The gradient projection algorithm with the metric
projection on the tangent plane}\label{GPA-Any}

The next version of the gradient projection algorithm uses the metric projection of the point
$x_k-tf'(x_k)$
on the tangent plane to the set $Q=\{ x\in\R^n\ |\ g(x)=0\}$, $g:\R^n\to\R$ is a $C^1$ function, at the point $x_k$. After this step we localize the next point $x_{k+1}$ on some segment and finding it
 by dividing the segment in half.

For a point $x$, $g(x)=0$, denote $T_x = \{ v\in\R^n\ |\ (g'(x),v)=
0
\}$ i.e. the tangent subspace to the surface $Q$ at the point $x\in Q$.

\begin{Lm}\label{LT2-1}
Assume that the function $f:\R^n\to\R$
is the Lipschitz continuous with constant $L>0$
and its gradient $f'$ is also the Lipschitz continuous with constant $L_1>0$.
Let $g:\R^n\to\R$ be a continuously differentiable function and
$Q = \{ x\in\R^n\ |\ g(x)=0\}$ be a surface without edge and
a proximally smooth set with constant $R>0$ with $g'(x)\ne 0$ for all $x\in Q$. Put
$$
t_0 = \frac{1}{L_1+\frac{2L}{R}}
$$
and fix $t\in (0, 2 t_0)$. Let $g(x_0)=0$,
$z=x_0-tP_{T_{x_0}}f'(x_0)$, $p=\frac{g'(x_0)}{\| g'(x_0)\|}$,
$$
\{x_1\}=\left[ z-p\left( R-\sqrt{R^2 - \| x_0 - z\|^2}\right),
z+p\left( R-\sqrt{R^2 - \| x_0 - z\|^2}\right)\right]\cap Q.
$$
Then
\begin{equation}\label{1step}
f(x_1)-f(x_0)\le -\| P_{T_{x_0}}f'(x_0)\|^2\left[ t-t^2\left(\frac{L_1}{2}+\frac{L}{R}\right)\right].
\end{equation}
The maximal value of the function $q(t) = t-t^2\left(\frac{L_1}{2}+\frac{L}{R}\right)$ is
$q_0 = q(t_0) = \frac12 t_0$, and $q(t) > 0$ for all $t \in (0, 2 t_0)$.
\end{Lm}

The Lemma implements following algorithm,
preferably with $t = t_0$.
\begin{framed}
\textbf{Gradient Projection on Tangent Hyperplane (GPA2)}
\vspace{3mm}

\textbf{Step 1.} Let $Q$ satisfy Lemma~\ref{LT2-1} condition.
Set $x_0 \in Q$, $0< t < 2 t_0$, and $k = 0$.

\textbf{Step 2.} Make a step and project onto tangent hyperplane:
$$z_k = x_k - t P_{T_{x_k}} f'(x_k)$$,

\textbf{Step 3.} Find intersection of a segment and the surface (i.e. by iterative bisection of the segment)
$$
\{x_{k+1}\} =
\begin{array}{l}
\left[
z_k - \frac{g'(x_k)}{\|g'(x_k)\|} \left(R - \sqrt{R^2 - \|x_k - z_k\|^2}\right) \right., \\ [3mm]
\left. z_k + \frac{g'(x_k)}{\|g'(x_k)\|} \left(R - \sqrt{R^2 - \|x_k - z_k\|^2}\right)
\right] \cap Q.
\end{array}
$$

\textbf{Step 4.}
Increase $k$ and continue to the Step 2.

\end{framed}

\begin{figure}[tbph]
    \centering
    \includegraphics[width=1\linewidth]{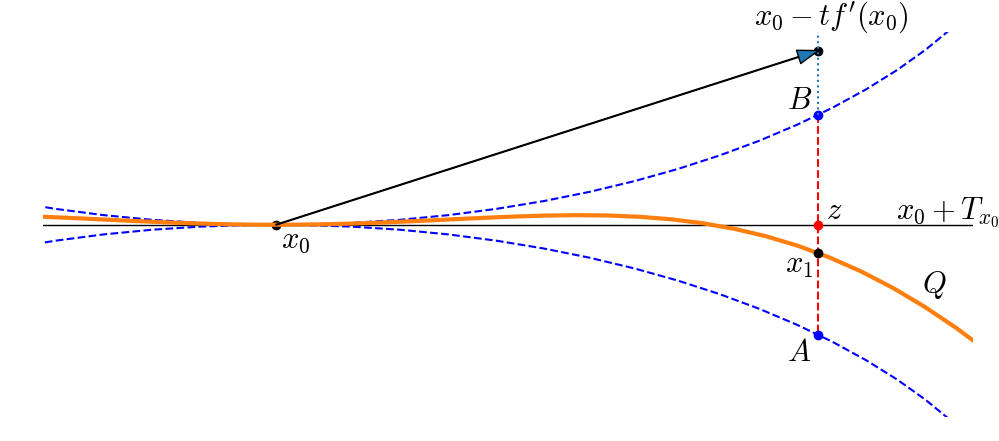}
    \caption{Lemma \ref{LT2-1}.}
    \label{fig3}
\end{figure}

\proof The maximality of $q_0$ is obvious. Let's prove (\ref{1step}).
By (\ref{BetweenS}) the segment
$$
[A,B]=\left[ z-p\left( R-\sqrt{R^2 - \| x_0 - z\|^2}\right),
z+p\left( R-\sqrt{R^2 - \| x_0 - z\|^2}\right)\right]
$$
has (unique) intersection $x_{1}$ with the set $Q$. The point $x_{1}$ can be found by
dividing the segment $[A,B]$ in half. See Figure \ref{fig3} for details.

We have
\begin{equation}\label{2-12}
f(x_1)-f(x_0)= f(x_1)-f(z)+f(z)-f(x_0),
\end{equation}
$\| x_0-z\| = t\| P_{T_{x_0}}f'(x_0)\|$,

\begin{equation}\label{2order}
\| x_1 - z\| \le R-\sqrt{R^2 - \| x_0 - z\|^2}\le \frac{\| x_0-z\|^2}{R} = \frac{t^2\| P_{T_{x_0}}f'(x_0)\|^2}{R}.
\end{equation}
\begin{equation}\label{2-1}
| f(x_1)-f(z)|\le L \| x_1-z\|\le \frac{t^2 L}{R}\| P_{T_{x_0}}f'(x_0)\|^2,
\end{equation}
and
$$
f(z)-f(x_0)\le (f'(x),z-x_0)+\frac{L_1}{2}\| z-x_0\|^2 = \| P_{T_{x_0}}f'(x_0)\|^2 \left( -t+\frac{L_1}{2}t^2\right).
$$
Substituting the last inequality and (\ref{2-1}) in Formula (\ref{2-12}) we get (\ref{1step}).\qed

Note that if the function $g$ in Lemma \ref{LT2-1} has the
Lipschitz continuous gradient with Lipschitz constant $L_{1}$ and
there exists $m>0$ with $\| g'(x)\|\ge m$ for all $x$, $g(x)=0$,
then the set $\{ x\ |\ g(x)=0\}$ is proximally smooth with
constant $R=m/L_{1}$ \cite[Proposition 4.15]{Vial}.

\begin{Th}\label{T2-1}
    Suppose that conditions of Lemma \ref{LT2-1} hold, $t\in (0, 2t_0)$, and the function $f$ satisfies
    the LPL condition with constant $\mu>0$ on the set $Q\cap \L_f (\beta)$.
    Then the GPA2 with initial condition
    $x_0\in Q\cap\L_f(\beta)$,
    converges with linear rate to the minimum point.
    \end{Th}
\proof Put $\varphi (x) = f(x)-f_0$, where $f_0 = \min_{Q}
f$. From the LPL condition for the function $f$ on the surface $Q$
$$
\| P_{T_x}f'(x)\|^2\ge \mu (f(x)-f_0)
$$
for all $x\in Q\cap\L_f(\beta)$ by Lemma \ref{LT2-1} we have
$$
f(x_{k+1})-f(x_k)\le -q(t)\mu (f(x_k)-f_0),
$$
$$
\varphi (x_{k+1})\le (1-q(t)\mu)\varphi (x_k).
$$
Now consider the rate of convergence with respect to the point. By (\ref{2order}) we get
$$
\| x_{k+1}-x_k\|^2 = \| x_k-z_k\|^2 + \| z_k-x_{k+1}\|^2 \le t^2\| P_{T_{x_k}}f'(x_k)\|^2
\left( 1+\frac{\| P_{T_{x_k}}f'(x_k)\|^2}{R^2}t^2\right).
$$
Using (\ref{1step}) for $x_0=x_k$ and $x_1=x_{k+1}$ we obtain that
$$
\| x_{k+1}-x_k\|^2 \le \left( 1+\frac{L^2}{R^2}t^2\right)\frac{t^2}{q(t)}(f(x_k)-f_0),
$$
$$
\| x_{k+1}-x_k\|^2 \le \left( 1+\frac{L^2}{R^2}t^2\right)\frac{t^2}{q(t)}\varphi (x_k) \le
\left( 1+\frac{L^2}{R^2}t^2\right)\frac{t^2}{q(t)} (1-q(t)\mu)^k\varphi (x_0).
$$
Due to inequalities for $t, q(t)$ this implies $\| x_{k+1}-x_k\|^2 \le cp^k, p<1$. The end of the proof is standard (compare the proof of Theorem 2).
\qed

\begin{Exa}\label{SCF}
Let $Q$ be a $C^{1}$ and proximally smooth with constant $R$
manifold without edge, $f:\R^{n}\to\R$ be a strongly convex function (with
constant of strong convexity $\varkappa>0$) with the Lipschitz
continuous gradient. Suppose that $f$ is the Lipschitz function
with constant $L>0$ on the level set $\L_{f}(\beta)$ and 
$\frac{L}{\varkappa}<R$.
Then the function $f$ satisfies the LPL
condition on the set $Q\cap\L_{f}(\beta)$. We shall give a sketch
of proof for this fact.
\end{Exa}

 By \cite[Lemma 2.1]{Balashov2017} the function $f$ has unique minimum $x_{*}\in
Q\cap\L_{f}(\beta)$. Fix a point $x_{0}\in
Q\cap\L_{f}(\beta)$ and put $z=x_{0}-\tau P_{T_{x_{0}}}f'(x_{0})$. Choose a positive number 
$\tau$  from the 
conditions of Lemma \ref{LT2-1},  $\tau$ less than $t$ from \cite[Formula
(8)]{Balashov2017} and $\tau\le \frac{R}{2L}$. Then by strong convexity of the function $f$ by \cite{Balashov2017}
we have linear rate of convergence for the GPA1 with step $\tau$. From Theorem 2.3 \cite{Balashov2017} for
$\x_{1}=P_{Q}(x_{0}-\tau f'(x_{0}))$ we get
$$
\| x_{0}-x_{*}\|\le \sum\limits_{k=0}^{\infty} q^k  \| x_0 - \x_{1} \|\le
\frac{\| x_{0}-\x_{1}\|}{1-q}
$$
where $q=q(\tau)\in(0,1)$ \cite[Formula (8)]{Balashov2017} and $q$
does not depend upon $x_{0}$.

Let $w=P_{x_{0}+T_{x_{0}}}\x_{1}$, $y=x_0-\tau f'(x_{0})$. Then from the definition of
$\x_{1}$ we have $\| y-x_0\|\ge \| y-\x_{1}\|$, thus we get the inequality $\| w-z\|\le \| x_{0}-z\|$. By
Formula (\ref{2order}) we obtain that $\| \x_{1}-w\|\le \frac1R \left( 2\|
x_{0}-z\|\right)^{2}$. Hence by the Pythagoras theorem
$$
\| \x_{1}-z\|^{2} = \| z-w\|^2 + \| w-\x_{1}\|^2\le\qquad\qquad\qquad\qquad\qquad\qquad\qquad\qquad\qquad\qquad
$$
$$
\qquad\qquad\qquad\qquad\qquad \le\| x_{0}-z\|^{2}+ \frac{1}{R^{2}} \left( 2\|
x_{0}-z\|\right)^{4}\le \| x_{0}-z\|^{2}\left( 1+
\frac{4\tau^{2}L^{2}}{R^{2}}\right),
$$
\begin{equation}\label{EBf}
\| x_{0}-x_{*}\|\le \frac{\| x_{0}-\x_{1}\|}{1-q}\le  \frac{\|
x_{0}-z\|+\| z-\x_{1}\|}{1-q}\le D \| x_{0}-z\| = D\tau\|
P_{T_{x_{0}}}f'(x_{0})\|,
\end{equation}
where $D = \left(1+\sqrt{1+ \frac{4\tau^{2}L^{2}}{R^{2}}}
\right)/(1-q)$.

Suppose that $f'(x_*)\ne 0$. Put $p=\frac{f'(x_{*})}{\| f'(x_{*})\|}$.  By the supporting
principle for proximally smooth sets
$$
\| x_{*}-x_{0}\pm Rp\|\ge R,
$$
$$
(p,x_{0}-x_{*})\le \frac{1}{2R}\| x_{0}-x_{*}\|^{2},\quad
(f'(x_{*}),x_{0}-x_{*})\le \frac{\| f'(x_{*})\|}{2R}\|
x_{0}-x_{*}\|^{2}.
$$
In the case $f'(x_*)= 0$ the last formula also takes place.
From (\ref{key}) and (\ref{EBf}) we obtain that
$$
f(x_{0})-f(x_{*})\le (f'(x_{*}),x_{0}-x_{*})+\frac{L_{1}}{2}\|
x_{0}-x_{*}\|^{2}\le \| x_{0}-x_{*}\|^{2}\left( \frac{\|
f'(x_{*})\|}{2R}+\frac{L_{1}}{2}\right)\le
$$
$$
\le\| x_{0}-x_{*}\|^{2}\left( \frac{
L}{2R}+\frac{L_{1}}{2}\right)\le D^{2}\tau^2\left( \frac{
L}{2R}+\frac{L_{1}}{2}\right)\| P_{T_{x_{0}}}f'(x_{0})\|^{2}.
$$

\subsection{The gradient method combined with the Newton method on the unit sphere}\label{GPA-NM}

Describe some symbiosis of the gradient projection algorithm and
the Newton method for finding a stationary point for the problem
$\min_Q f$. We shall assume that $f\in C^2$.

Consider again the problem with $Q=S_{1}$. Define $S_{1}$ with the
help of the function $g(x) = \frac12 (\| x\|^2-1)=0$. For any $\|
x\|=1$ define the number $\lambda \in\R$ as a solution of the extremal problem
$$
\| f'(x)+\lambda g'(x)\|^2\to \min\limits_{\lambda},
$$
thus $\lambda = -(g'(x)^Tg'(x))^{-1}g'(x)^Tf'(x) = - (x,f'(x))$.
Denote by $z$ the variable $(x,\lambda)\in\R^{n+1}$,
$$
F(z) = \left[\begin{array}{c} f'(x) + \lambda x \\
g(x) \end{array}\right],
\quad
F'(z) = \left[\begin{array}{cc} f''(x) + \lambda I & x\\
x^{T} & 0 \end{array} \right].
$$

Fix $z=(x,\ll)$, where $ x\in S_{1}$, $\ll = -(x,f'(x))$. Define
also
$$
\sigma_1 = \arg\min\limits_{\sigma\in \sigma (F'(z)) }|\sigma|
$$
the minimal by absolute value element of spectrum for the matrix $F'(z)$.

Suppose that $\sigma_1>0$ and $F'$ is the Lipschitz continuous
with Lipschitz constant $L_1>0$ on the set $B_{r}(z)$, where
$r=\frac{2}{\sigma_{1}}\| F(z)\|$. Then in the case
\begin{equation}\label{MNew}
\frac{L_1}{\sigma_1^2}\sqrt{\| f'(x)+ \ll x\|^2 +g^2(x) }=
\frac{L_1}{\sigma_1^2}\sqrt{\| f'(x)- (x,f'(x))x\|^2 +g^2(x)
}<\frac14
\end{equation}
the modified Newton method starting from the point
$z=(x,-(x,f'(x)))$ converges with super-linear rate \cite[Chapter
X, \S 4, Theorem 1]{KF}.

Note that $\|P_{T_x}f'(x)\| =  \|f'(x)-(x,f'(x))x\| = \|
(I-xx^{T})f'(x)\|$.

\begin{framed}
    \textbf{Gradient Projection --- Newton Method (GPA3)}

    \textbf{Step 1.}
    Take $x_0 \in {S_{1}}$, and put $C = \frac{\sigma_1^2}{4 L_1}$, $k = 0$.

    \textbf{Step 2.} (GPA2 phase) While $\|P_{T_{x_k}}f'(x_k)\| \ge C$,
    perform Steps~2-4 of GPA2, increasing $k$.
    If $\|P_{T_{x_k}}f'(x_k)\| < C$, proceed to Step~3.

    \textbf{Step 3.} Put $x_{0}=x_{k}$, $\lambda_0 = -(x_k,f'(x_k))$,
    $k=0$.

    \textbf{Step 4.} (Newton phase) do Newton steps for equation $F(z) = 0$, increasing $k$:
    $$
    \begin{array}{rl}
    \begin{bmatrix}
    x_{k+1} \\
    \lambda_{k+1}
    \end{bmatrix}
    & =
    \begin{bmatrix}
    x_{k} \\
    \lambda_{k}
    \end{bmatrix}
    - F'\Big(\begin{bmatrix}
    x_{0} \\
    \lambda_{0}
    \end{bmatrix}
    \Big)^{-1}
    F\Big(\begin{bmatrix}
    x_{k} \\
    \lambda_{k}
    \end{bmatrix}
    \Big) = \\
    & =
    \begin{bmatrix}
    x_{k} \\
    \lambda_{k}
    \end{bmatrix}
    -
    \begin{bmatrix}
    f''(x_0) + \lambda_0 I & x_0\\
    x_0^T & 0
    \end{bmatrix}^{-1}
    \begin{bmatrix}
    f'(x_k) + \lambda_k x_k \\
    g(x_k)
    \end{bmatrix}
    \end{array}
    $$

\end{framed}
Conditions of Lemma \ref{LT2-1} are satisfied at Step~2.
Thus we can do steps (\ref{1step}) of the GPA2
and decrease the function:
$$
f(x_{k+1})-f(x_k) \le - C^2 q(t).
$$
Put $\Delta f = \max_{S_{1}} f - \min_{S_{1}} f$. It's easy to
see that by the inequality $\|P_{T_{x_k}}f'(x_k)\| \geq C$ we'll
switch to the Newton method at Step~4 after no more than
$\frac{\Delta f}{C^2q(t)}$ steps of the gradient projection
algorithm. In the case when condition (\ref{MNew}) is valid at
the point $z=(x_{k},-(x_{k},f'(x_{k})))$ the modified Newton
method starting from $z$ converges with super-linear rate.

\subsection{Quadratic form}\label{QF}

Consider homogeneous quadratic function $f(x)=(Ax,x)$ with symmetric
real $n\times n$ matrix $A$. Denote by $\lambda_k$ eigenvalues of
$A$, $e_k$ --- corresponding eigenvectors $A e_k = \lambda_k e_k$
and suppose that $\lambda_1 < \lambda_2 \le \dots \le \lambda_n$.
Then $f_0=\lambda_1, L_{1}= 2\lambda_n$ and two global
minimizers are $\pm e_1$. All other eigenvectors are stationary
points, but not local minimums. Thus the problem
(\ref{sphere}) is equivalent to the problem of finding the minimal eigenvalue and
eigenvector, and algorithm (\ref{alg-with-L}) has the form
\begin{equation}\label{eig}
x_{k+1} = \frac{ x_k -(2/L_{1})Ax_k}{\| x_k -(2/L_{1}) Ax_k\|}.
\end{equation}

Probably, the first gradient-like algorithm for eigenvalue
problem has been proposed by Kantorovich \cite[Section
3.4]{Kantorovich}. He converted eigenvalue problem to
unconstrained minimization of Rayleigh quotient
$R(x)=\frac{(Ax,x)}{(x,x)}$ and obtained the algorithm $x_{k+1} = \frac{
x_k -t_k Ax_k}{\| x_k -t_k Ax_k\|}$ where $t_k$ was taken from 1D
minimization of $R(x_{k+1})$. One can see that this method has
the same form as (\ref{eig}), but has more complicated
step-size rule. Kantorovich proved linear convergence of the
algorithm.

We analyse iterative process (\ref{eig}) by use of the above presented
results. For $x_0 = \sum_{i=1}^{n} \alpha_i e_i$ we shall prove
convergence to $+e_1$ or $-e_1$ depending on the sign of $\alpha_1$.
Suppose that $\alpha_1>0$ (if $\alpha_1=0$, there is no convergence
to global minimum). It is obvious that if $\alpha_1>0$, the same
is true for decompositions of all iterates $x_k$, thus all of
them remain in the open half-sphere $H=\{ x\ |\ (x,e_1)>0,
||x||=1\}$. We also introduce $H_{\tau}=\{ x\ |\  (x, e_1) \geq
\tau, \; ||x|| = 1\},$ for $0 < \tau <1$.

\begin{Lm}\label{Axx} Fix $\tau\in (0,1)$. Suppose that
    $\lambda_1 < \lambda_2 \leq ... \leq \lambda_n$ are eigenvalues
    of $A$.
Then the quadratic function $f(x) = (A x, x)$ satisfies the LPL condition
on the set $H_{\tau}$ with $\mu = 4 \tau^2(\lambda_2 - \lambda_1)$,
\end{Lm}

\proof Express any point $x = \sum_{i=1}^{n} \alpha_i e_i \in
H_\tau$ through the residual vector $\delta$:
$\delta = x -e_1$.
From $\tau\in (0,1)$ we have $0<\|\delta\| < \sqrt{2}$. Put $C=\| \delta\|^{2}-\frac{\|
\delta\|^{4}}{4}>0$. Using notation $y$ for the unit vector $y =
\frac{1}{\sqrt{C}}[\alpha_2, \alpha_3, ..., \alpha_n]^T \in
\mathbb{R}^{n-1}$ (see (\ref{sum_of_alpha_2_n})), and $B$ for the
diagonal matrix with strictly positive diagonal elements
$$
B = \begin{bmatrix}
\lambda_2 - \lambda_1 & 0 & \cdots & 0 \\
0 & \lambda_3 - \lambda_1 & \cdots & 0 \\
\vdots & \vdots & \ddots & \vdots \\
0 & 0 & \cdots & \lambda_n - \lambda_1 \\
\end{bmatrix}
$$
we have the next obvious equalities below

\begin{gather}
\nonumber 1 = \|x\|^2 = \sum_{i=1}^{n} \alpha_i^2, \qquad
\label{delta_norm-and-alpha_1} \|\delta\|^2 = (1 - \alpha_1)^2 +
\sum_{i=2}^{n} \alpha_i^2 = 2 - 2 \alpha_1,\\
\label{sum_of_alpha_2_n} \sum_{i=2}^{n} \alpha_i^2 = 1 -
\alpha_1^2 = \|\delta\|^2 - \frac{\|\delta\|^4}{4} \doteq
C,\qquad (A - \lambda_1 I) \delta = \sum_{i=2}^{n} (\lambda_i -
\lambda_1) \alpha_i e_i,\\
\nonumber \label{squared-norm}
\|(A - \lambda_1 I) \delta\|^2 = \sum_{i=2}^{n} (\lambda_i - \lambda_1)^2 \alpha_i^2 = C \|B y\|^2,\\
\label{delta_f}
((A - \lambda_1 I) \delta, \delta) = \sum_{i=2}^{n} (\lambda_i - \lambda_1) \alpha_i^2 = C (B y, y), \\
\label{f_x-f*}
f(x) - f_0 = (Ax, x) - \lambda_1 = \sum_{i=1}^{n}(\lambda_i - \lambda_1)\alpha_i^2 =
\sum_{i=2}^{n}(\lambda_i - \lambda_1)\alpha_i^2 = C (By, y).
\end{gather}
It is clear that $(By, y) \geq \lambda_2 - \lambda_1 > 0$.

From the equality $f'(x) = 2Ax = 2A(e_1 + \delta)$ we get
\begin{gather*}
\|(I - x x^T) f'(x)\|^2 = \|(I - (e_1 + \delta) (e_1 + \delta)^T) (2A (e_1 + \delta))\|^2 = \\
= 4 \Big\|(A - \lambda_1 I)\delta
+ (e_1 + \delta) (A \delta - \lambda_1 \delta, \delta) \Big\|^2
= \\
= 4\|(A - \lambda_1 I) \delta\|^2 - 4 ( (A - \lambda_1 I)\delta, \delta )^2 = \\
= 4C \|B y\|^2 - 4 C^2 (By, y)^2 = 4 C (\|B y\|^2 - C (By, y)^2) = \\
= 4 C (\|B y\|^2 - (By, y)^2) + 4 C (1 - C) (By, y)^2 \geq \\
\geq 4 (1 - C) (By, y) \big[C (By, y)\big] = 4 (1 - C) (By, y) (f(x) - f_0).
\end{gather*}
By $4(1-C) = (2 - \|\delta\|^2)^2$ the latter expression has the form
$$
\|(I - x x^T) f'(x)\|^2 \geq (2 - \|\delta\|^2)^2 (\lambda_2 - \lambda_1) (f(x) - f(x^*)).
$$
for all $x\in H_{\tau}$.
The inequality $2 - \|\delta\|^2 = 2 - \|x - e_1\|^2 = 2 (x, e_1)
\geq 2 \tau$ holds for $x \in H_\tau$, and thus the LPL condition
takes place with $\mu = 4 \tau^2(\lambda_2 - \lambda_1)$.\qed

Thus for any $(x_0, e_1) \ge \tau$ we get $(Ax_k, x_k) -
\lambda_1 \le ((A x_0, x_0) - \lambda_1) q^k, q = 1 - \tau^2
\frac{\lambda_2 - \lambda_1}{\lambda_n - \lambda_1}$ (denominator
$\|L x_k - f'(x_k)\|$ from \eqref{formula-for-q} is bounded from
above by $2(\lambda_n - \lambda_1)$ for the quadratic function),
while asymptotically $|(x_k, e_1)| \rightarrow 1$, and
$(Ax_k,x_k)-\lambda_1=O (q_1^k),\; q_1 = \frac{\lambda_n -
\lambda_2}{\lambda_n - \lambda_1}$.

Condition $\lambda_1<\lambda_{2}\dots$ can be weakened.
If $\lambda_1=\dots=\lambda_k<\lambda_{k+1}\le \dots\le \lambda_n$, then the function
$f(x)=(Ax,x)$ also satisfies the LPL condition on the set
$H_{\tau} = \{ x\in S_1\ |\  \sum_{i=1}^{k}x_i^2\ge \tau^2\}$
 of the unit sphere $S_1$
 for any $\tau\in (0,1)$
 in the basis from eigenvectors.

\def\x {\overline{x}}

Indeed, put $z=(x_1,\dots,x_k)$, $\x = (x_{k+1},\dots,x_n)$, $x=(z,\x)$.
Then
$$
f(x) = \sum\limits_{i=1}^{n}\lambda_ix_i^2 = \lambda_1 \|z\|^2+
\sum\limits_{i=k+1}^{n}\lambda_ix_i^2 = h(\| z\|,\x).
$$
For the function $h(\| z\|,\x)$ we have the LPL condition on the set
$\{ \| z\|> \tau\}$ with constant $\mu = 4\tau^2(\lambda_{k+1}-\lambda_1)$, i.e.
$$
4\left( \lambda_1^2\|z\|^2+\sum\limits_{i=k+1}^{n}\lambda_ix_i^2\right)\ge \mu (h(z,\x)-f_0),
$$
or equivalently
$$
4\left( \sum\limits_{i=1}^{n}\lambda_ix_i^2\right)\ge\mu (f(x)-f_0).
$$

\section{The Frank-Wolfe method}\label{FW-sec}

\emph{ The Frank-Wolfe method} (also known as \emph{the conditional gradient method}) has been proposed for minimization of a convex quadratic
function on a convex set \cite{FW} and later was extended for general convex objectives, see e.g. \cite{Polyak_book} and recent survey \cite{FW_1}.
The idea of the method for problem (\ref{sphere}) is to solve (on each step) the
auxiliary problem
$$
z_k\in\mbox{\rm Arg}\min\limits_{z\in Q}(f'(x_k),z),
$$
find
$t\in [0,1]$ that minimizes $f(x_k+t(z_k-x_k))$ and take the next point
$x_{k+1}=x_k+t(z_k-x_k)$.
The method requires minimization of a linear function on the admissible set at each iteration. There are also some extensions of the
method for nonconvex objective functions and for matrix optimization \cite{Lacoste-2016,Weber-Sra-}.

\subsection{Minimization of an approximately linear function}\label{FW-ssec1}

However our problem (\ref{sphere})
deals with non-convex admissible set $Q$. We consider a special version of the FW method for our problem as a limiting version of the
gradient projection method (\ref{alg-with-L}). Indeed suppose that the function $f(x)$ is approximately linear (see (\ref{appr.lin.}) below) on the set $Q=S_1$, i.e., informally, constant $L_1$ is
small enough in comparison with other parameters.
 For this extreme case method (\ref{alg-with-L}) turns into the next iteration process
\begin{equation} \label{full-fw-on-sphere}
x_{k+1} = -\frac{f'(x_k)}{\|f'(x_k)\|}.
\end{equation}
This is exactly the FW method with $t=1$: we take linearized
function $f(x_k)+(f'(x_k), x-x_k)$, find its minimum on $Q$ and
proceed to the minimum point. Notice that in standard versions
of the FW method we make a step in the direction of the
minimizer;  this full-step version diverges in general case.

\begin{framed}
\textbf{Full-step Frank-Wolfe method (FFW)}

\textbf{Step 1.} Take $x_0$ and set $k = 0$.

\textbf{Step 2.} Solve auxiliary problem
$$
z_k \in \mbox{\rm Arg}\min\limits_{z\in Q}(f'(x_k),z),
$$

\textbf{Step 3.} Update
$x_{k+1} = z_k.$

\textbf{Step 4.} Increase $k$ and go to Step 2.
\end{framed}

To get the rigorous validation of method (\ref{full-fw-on-sphere}) we need
specification of the above mentioned approach.
A function $f(x)$ defined on the ball $B_1(0)$ with $L_1$-Lipschitz gradient (\ref{Lip}) is called \emph{approximately linear} on $S_1$ if
\begin{equation}\label{appr.lin.}
||f'(0)|| > 2 L_1.
\end{equation}

\begin{Th} \label{FW-global}
Suppose that (\ref{appr.lin.}) holds. There are just two
stationary points in problem (\ref{sphere}) $x_{*} =
\arg\min_{S_{1}}f(x),$ $ x^{*}=\arg\max_{S_{1}} f(x)$, and
FFW method \eqref{full-fw-on-sphere}
converges to $x_{*}$
for arbitrary $x_0\in S_1$ with linear rate
\begin{equation}\label{4}
||x_k-x_{*}||\le q^k ||x_0-x_{*}||, \quad q=
\frac{L_{1}}{\|f'(0)\| - L_{1}}\in (0,1).
\end{equation}
\end{Th}
Theorem \ref{FW-global} is close enough to Theorem~4.3 from
\cite{Polyak2001} where minimization on $B_1$ instead of $S_1$ has been considered. But indeed the solutions of these two problems coincide under condition \eqref{4}.
The proof of Theorem \ref{FW-global} follows from the following
 fact regarding strongly convex sets of radius $r$ and
functions with the Lipschitz continuous gradient
\cite{Balashov2015}.

Suppose that $B\subset\R^n$ is a strongly convex set of radius $r$
and a function $f:\R^n\to\R$ has the Lipschitz continuous
gradient with constant $L_1$, $m=\inf\limits_{x\in \d Q}\|f'(x)\|$ and
$\frac{m}{L_1}>r$. Then the iteration process $x_0\in\d B$,
$$
x_{k+1}=\arg\min\limits_{x\in B}(f'(x_k),x),\qquad k=0,1,\dots
$$
converges to the unique solution $x_*\in \d B$ of the problem
$\min_B f$ with linear rate:
$$
\| x_k -x_*\|\le \left( \frac{rL_1}{m}\right)^k\| x_0 - x_*\|.
$$

We shall further consider a closed surface $Q$ in $\R^{n}$ which is the
boundary of some strongly convex set of radius $r$, i.e. $Q=\d B$.
It is worth to admit that $Q$ is not necessary smooth.

\subsection{Another gradient domination condition}\label{FW-ssec2}

We introduce a sort of the gradient domination condition, formulated at a stationary point $x_*$ of (1).
Assume that the next inequality
\begin{equation}\label{almost-linear-2}
||f'(x_{*})|| > L_1
\end{equation}
holds. This condition reminds \textit{sharp minimum} condition \cite{Polyak_book}.

\begin{Th}\label{NesCond}
    Let $B\subset \R^{n}$ be a strongly convex set of radius $r$,
    $m>1$, $x_{*}\in\d B$. Suppose that $f:\R^{n}\to\R$ is a function
    with the Lipschitz continuous gradient with constant $L_{1}$. If
    $f'(x_{*})\in -N(B,x^{*})$ and $\| f'(x_{*})\|\ge mrL_{1}$ then
    $x_{*}=\arg\min\limits_{a\in B}f(a)$, i.e. $x_{*}$ is the strict
    global minimum of the function $f$ on the set $B$, and hence on the boundary $Q=\d B$.
\end{Th}

\proof Put $p=\frac{f'(x_{*})}{\| f'(x_{*})\|}$. Then by the
supporting principle for strongly convex sets
\begin{equation}\label{0*}
B\subset B_{r}(x_{*}+rp).
\end{equation}
Fix a number $\ep>0$ and a unit vector $q$ such that
\begin{equation}\label{1*}
x=x_{*}+\ep q\in B_{r}(x_{*}+rp).
\end{equation}
We claim that $f(x)>f(x_{*})$.

By Formula (\ref{key})
\begin{equation}\label{2*}
f(x)\ge f(x_{*})+(f'(x_{*}),x-x_{*})-\frac{L_{1}}{2}\|
x-x_{*}\|^{2}.
\end{equation}
From (\ref{1*}) we get $\| rp-\ep q\|\le r$. Hence $-2r\ep
(p,q)+\ep^{2}\le 0$ or $\frac{\ep}{2}\le r(p,q)$. By inequality $m>1$ we obtain that
\begin{equation}\label{3*}
\frac{\ep}{2}\le m r(p,q).
\end{equation}
Further
$$
(f'(x_{*}),x-x_{*}) = (f'(x_{*}),\ep q) \ge \ep m rL_{1}(p,q)
$$
and $\| x-x_{*}\|^2=\ep^{2}$. The last two formulae and (\ref{2*})
gives the next estimate
$$
f(x)\ge f(x_{*})+ \ep m rL_{1}(p,q) -
\frac{L_{1}}{2}\ep^{2}=f(x_{*}) + \ep L_{1}\left(
mr(p,q)-\frac{\ep}{2}\right),
$$
and taking in mind (\ref{3*}) we have $f(x)>f(x_{*})$. By
inclusions (\ref{0*}), (\ref{1*}) $x_{*}=\arg\min\limits_{a\in B}f(a)$.\qed

\begin{Corol}\label{m=1}
    If $m=1$ and other assumptions of Theorem \ref{NesCond} hold then
    $x^{*}$ is also a global minimum of $f$ on the set $B$, but this
    minimum is not necessary strict.
\end{Corol}

The particular example is given by the set $B=B_{1}(0)$ and the
function $f(x)=-\frac12 \|x\|^{2}$.

\begin{Corol}\label{m<1}
    If in Theorem \ref{NesCond} $m\in (0,1)$
    and other assumptions of
    Theorem \ref{NesCond} hold then $x_{*}$ is a stationary point of
    the function $f$ on the set $B$ but not necessarily the minimum point.
\end{Corol}

\begin{Exa}\label{MinStat}
Consider an example in $\R^{2}$. Fix $r>1$, let $B=\{ (x,y)\in
R^{2}\ |\ x^{2}+(y+r)^{2}\le r^{2}\}=B_{r}((0,-r))$.

Define the function $\psi$
$$
\psi (x) = \left\{\begin{array}{c} -\frac12 x^{2},\quad x\le 0,\\
0,\quad x> 0,\\
\end{array}\right.
\qquad
\psi' (x) = \left\{\begin{array}{c} -x,\quad x\le 0,\\
0,\quad x> 0\\
\end{array}\right.
$$
and $f(x,y) = \psi (x)-y$ for all $(x,y)\in\R^{2}$.
\end{Exa}
Consider the
problem $\min_{B}f$. We have $f'(x,y) = (\psi'(x),-1)$ and
thus $f'$ is the Lipschitz continuous with constant $L_{1}=1$.

The point $a=(0,0)\in \d B$ is a stationary point: $f'(a)
=(0,-1)\in -N(B,a)$, $\| f'(a)\| = 1 = mrL_{1} = mr$, i.e.
$m=\frac1r<1$ and $f(a)=0$. But  $a$ is not a local
minimum. The solution of the problem is $a_{0} =
(-\sqrt{r^{2}-1},1-r)$ with $f(a_{0}) = -\frac{(r-1)^{2}}{2}$.\qed

Define the function $h(\theta) = 2\sin\left(\frac12\arcsin \theta\right)$ for $\theta\in [0,1]$. The function $h$ is convex, monotonically increasing and $\frac{1}{\theta}h(\theta)>1$
for $\theta >0$. By convexity of $h$ we have for any $t\in (0,\theta)$ that
$h(t)\le \frac{1}{\theta}h(\theta) t$. Note that $1\le \frac{1}{\theta}h(\theta)
\le\sqrt{2}$ if $\theta\in [0,1]$.

For any real number $m>1$ define $\theta_m$ as follows
$$
\theta_m = \left\{\begin{array}{lc}
\theta, \quad \mbox{where}\ \frac{1}{\theta}h(\theta)=m\ \mbox{and}\ m\in (1,\sqrt{2}],\\
1, \ m>\sqrt{2}.
\end{array}\right.
$$

\begin{Th}\label{FWx*1}
    Let $B\subset \R^{n}$ be a strongly convex set of radius $r$,
    $m>1$, $x_{*}\in\d B$. Suppose that $f:\R^{n}\to\R$ is a function
    with the Lipschitz continuous gradient with constant $L_{1}$,
    $f'(x_{*})\in -N(B,x_{*})$ and $\| f'(x_{*})\|\ge mrL_{1}$. Fix a
    point $x_{0}\in \d B$, $\| x_{0}-x_{*}\|<
    \theta mr$ and $\theta\in (0,\theta_m)$. Then the
    iterations
    $$
    x_{k+1}=\arg\max\limits_{x\in \d B}(-f'(x_{k}),x),\quad
    k=0,1,\dots
    $$
    of FFW method converge to the point $x_{*}$ with linear rate:
    $$
    \| x_{k+1}-x_{*}\| \le \frac{1}{\theta}\frac{h(\theta)}{m}\| x_{k}-x_{*}\|
    $$
    for all $k$.
\end{Th}

Note that by Theorem \ref{NesCond} $x_{*}$ is the strict global
minimum of $f$ on the set $B$. Also
$\frac{1}{\theta}\frac{h(\theta)}{m}<1$ due to choice of $\theta$.

\proof Put $t_{k}=\| x_{k}-x_{*}\|$.

From the inequality $L_{1}t_{k}<\|
f'(x_{*})\|$ (which follows by induction) and the inclusion $f'(x_k)\in B_{L_1t_k}(f'(x_*))$ the sine of the angle $\varphi_k$
between $f'(x_*)$ and $f'(x_k)$ is estimated as follows
$$
\sin\varphi_k \le \frac{L_1t_k}{\| f'(x_*)\|}\le \frac{t_k}{mr}.
$$
Put $p_k=-\frac{ f'(x_k)}{\| f'(x_k)\|}$, $p_*=-\frac{ f'(x_*)}{\| f'(x_*)\|}$.
From the triangle $0p_kp_*$ we have $\| p_k-p_*\|=2\sin\frac12\varphi_k$.
By strong convexity of the set $B$ and inequality $\frac{t_k}{mr}<\theta$ we obtain that
$$
t_{k+1}=\| x_{k+1}-x_*\|\le r\| p_k-p_*\|=2r\sin\frac12\varphi_k \le
2r\sin\left(\frac12\arcsin\frac{t_k}{mr}\right)\le \frac{1}{\theta}\frac{h(\theta)}{m}t_k.
$$
\qed

\begin{Th}\label{FWx01}
    Let $B\subset \R^{n}$ be a strongly convex set of radius $r$,
    $m>1$, $\beta\in\R$. Suppose that
    $f:\R^{n}\to\R$ is a function with the Lipschitz continuous
    gradient with constant $L_{1}$ and for any point $x\in \d B\cap\L_f(\beta)$
we have the inequality  $\| f'(x)\|\ge mrL_{1}$.
    Then for any choice of the initial point $x_0\in \d B\cap\L_f(\beta)$
    the iterations
    $$
    x_{k+1}=\arg\max\limits_{x\in \d B}(-f'(x_{k}),x),\quad
    k=0,1,\dots
    $$
    converge to the global strict minimum $x_{*}$ with linear rate:
    $$
    \| x_{k+2}-x_{k+1}\| \le \frac{1}{m}\| x_{k+1}-x_k\|.
    $$
    for all $k$.
\end{Th}

\proof Put $t_k=\| x_{k+1}-x_k\|$ and $p_k=-\frac{f'(x_k)}{\| f'(x_k)\|}$.

Suppose that $f(x_i)\le \beta$ for all $0\le i\le k$.

Prove that $f(x_{k+1})\le \beta$.
$$
f(x_{k+1})-f(x_k)\le (f'(x_k),x_{k+1}-x_k)+\frac{L_1}{2}\| x_{k+1}-x_k\|^2.
$$
From the supporting principle for strongly convex sets
$$
\| x_{k+1}-rp_k - x_k\| <r,
$$
or
$$
\| x_{k+1}-x_k\|^2\le 2r(p_k, x_{k+1}-x_k)\le -\frac{2}{mL_1}(f'(x_k), x_{k+1}-x_k).
$$
Hence we obtain the next estimate
$$
f(x_{k+1})-f(x_k)\le -\left(m-1 \right)\frac{L_1}{2}\| x_{k+1}-x_k\|^2,
$$
and $f(x_{k+1})\le f(x_k)\le \beta$.

For any unit vectors $p,q$ and numbers $\lambda,\mu\ge 1$ we have $\| p-q\|\le
\|\lambda p - \mu q \|$. Using the last inequality and strong convexity of the set $B$ we get
$$
t_{k+1}=\| x_{k+2}-x_{k+1}\|\le r\left\|\frac{f'(x_{k+1})}{mrL_1}- \frac{f'(x_{k})}{mrL_1}\right\|
\le \frac{1}{m}t_k.
$$
Thus the sequence $\{ x_k\}$ converges, $x_k\to x_*\in\d B$.
Passing to the limit as $k\to\infty$ in the inclusion
$-f'(x_{k})\in N(B,x_{k+1})$ and using upper semicontinuity of the
normal cone we have $-f'(x_{*})\in N(B,x_{*})$, $\| f'(x_*)\|\ge
mrL_1$. By Theorem \ref{NesCond} the point $x_*$ is the strict global
minimum. \qed

Consider examples that show importance of the condition $m>1$ in both
Theorems \ref{FWx*1} and \ref{FWx01}.

\begin{figure}[tbph]
    \centering
    \includegraphics[width=0.5\linewidth]{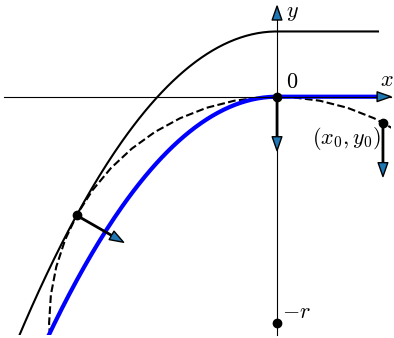}
    \caption{Stationary point $(0,0)$ is not a solution.}
    \label{fig2}
\end{figure}

\begin{Exa}\label{E1} \rm
Suppose that $f(x,y)$ and $B$ are the function and the set from
Example \ref{MinStat}. Notice that in this case $m<1$. Take a
starting point for the FFW algorithm $(x_{0},y_{0})\in\d B$ with
$x_{0}>0$, $y_{0}\in (-1,0)$. Then $f'(x_{0},y_{0})=(0,-1)$ and
one step of the FFW method leads us to the stationary point
$(0,0)$. But extremum is the point $(-\sqrt{r^{2}-1},1-r)$, see Figure \ref{fig2}.
\end{Exa}

\begin{Exa}\label{E2} \rm
Consider the set $Q=\{ (x,y)\in\R^{2}\ |\
x^{2}+(y-\frac12)^{2}=\frac14\}$. Let $f(x,y) =
y-x^{2}$.
For the function $f$ we have $L_1=2$, $1=\| f'(0,0)\|=mrL_1 =m\frac12\cdot 2=m$, i.e. $m=1$,
$$
(0,0) = \arg\min\limits_{(x,y)\in Q}f(x,y)
$$
and $f_{0} = f(0,0)=0$. Put $t>0$ and $a_{0}=(x_{0},y_{0})\in
Q$, $x_{0}>0$, $y_{0}\in (0,\frac12)$ with $f(x_{0},y_{0}) = y_{0}-x_{0}^{2}=t$. As we've seen at
subsection \ref{11} the
angle between tangent lines to the circle $Q$ and curve
$f(x,y)=y-x^{2}=t$ at the point $a_{0}$ asymptotically equals
$x_{0}^{3}$ when $x_{0}\to 0$.
Starting the FFW algorithm from the point
$a_{0}$ we obtain the next point $a_{1}(x_{1},y_{1})$. We have
$x_0^3\asymp\| a_{0}-a_{1}\|\asymp | x_{0}-x_{1}|$ when $x_{0}\to 0$. The
last means that
$$
\exists \delta>0 \ \exists C>0\ \forall x_0\in (0,\delta)\ |
x_{0}-x_{1}|<Cx_{0}^{3}.
$$
Thus for any sufficiently small $x_{0}>0$ we get $x_{1}\ge
x_{0}-Cx_{0}^{3}$. There is no linear rate of convergence.
\end{Exa}

%
%

\section{Appendix}

\begin{Prop}\label{NecMin}
    Suppose that the set $Q\subset\R^{n}$ is proximally smooth with
    constant $R>0$, the function $f$ has the Lipschitz continuous
    gradient $f'$ with constant $L_{1}>0$ and the point $x_{0}$ is a local
    minimum in the problem $\min_{Q}f$.
    Then $-f'(x_{0})\in N(Q,x_{0})$.
\end{Prop}

\proof
If $f'(x_0)=0$ then $-f'(x_0)=0\in N(Q,x_0)$ by the definition of normal cone.

Assume that $f'(x_{0})\ne 0$.
Prove the proposition by contradiction.
Put $p=-\frac{f'(x_{0})}{\|f'(x_{0})\|}$.
Suppose that $p\notin N(Q,x_{0})$. Then by the
supporting principle for a proximally smooth set we have
$$
\exists x_{1}\in Q\cap \i B_{R}(x_{0}+Rp).
$$

\begin{figure}[tbph]
    \centering
    \includegraphics[width=0.6\linewidth]{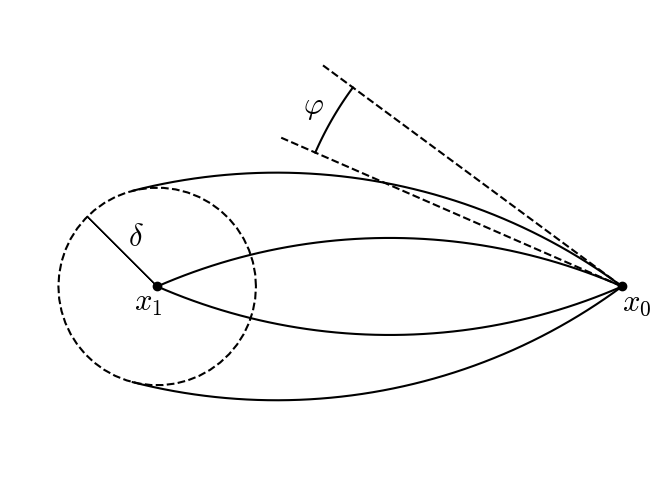}
    \caption{Angle $\varphi$.}
    \label{fig1}
\end{figure}

For a set $B$ define by $D_{R}B$ the intersection of all closed balls of radius
$R$ each of which contains the set $B$. From \cite[Lemmata 4.13,
4.16]{BI} there exists a continuous curve $\Gamma\subset Q$ with endpoints $x_0$
and $x_1$ such that $\Gamma\subset
D_{R}\{x_{1},x_{0}\}\subset \{ x_{0}\}\cup\i B_{R}(x_{0}+Rp)$. By the inclusion $x_{1}\in \i
B_{R}(x_{0}+Rp)$ there exists $\delta>0$ with
$B_{\delta}(x_{1})\subset \i B_{R}(x_{0}+Rp)$,
$$
D_{R}\{x_{1},x_{0}\}\subset D_{R}(B_{\delta}(x_{1})\cup
\{x_{0}\}).
$$
Let $L$ be a 2-dimensional plane, $\{ x_{0},x_{1}\}\subset L$.
Choose a point $x\in\Gamma$, $\| x-x_{0}\|=\ep$, $x=x_{0}+\ep q$,
$\| q\|=1$. The angle $\varphi$ between arcs $L\cap \d
D_{R}\{x_{1},x_{0}\}$ and $L\cap \d D_{R}(B_{\delta}(x_{1})\cup
\{x_{0}\})$ at the point $x_{0}$ (see Figure \ref{fig1})
is strictly positive
and hence the angle between $p$ and $q$ is less than $\frac{\pi}{2}-\varphi\in
(0,\frac{\pi}{2})$. Thus there exists $C= \cos\left( \frac{\pi}{2}-\varphi\right)>0$ (that does not depend on vector $q$)
with $(p,q)>C$.

By Formula (\ref{key})
$$
f(x) - f(x_{0})\le (f'(x_{0}),x-x_{0})+\frac{L_{1}}{2}\|
x-x_{0}\|^{2}
$$
and
$$
(f'(x_{0}),x-x_{0}) = (f'(x_{0}),\ep q) = - \| f'(x_{0})\|\ep
(p,q).
$$
Hence
$$
f(x) - f(x_{0})\le - \| f'(x_{0})\|\ep (p,q) +
\frac{L_{1}}{2}\ep^{2} \le - \| f'(x_{0})\|\ep C +
\frac{L_{1}}{2}\ep^{2}<0
$$
for sufficiently small $\ep>0$. A contradiction.\qed

\section*{Acknowledgements}

The work was supported by Russian Science Foundation (Project 16-11-10015).


\end{document}